\documentclass{article}

\usepackage{hyperref}
\usepackage{a4wide}
\usepackage{amsfonts,amssymb,amsmath,amsthm}
\usepackage{graphicx}
\usepackage[latin1]{inputenc}
\usepackage{color}
\usepackage{epstopdf}
\usepackage{subfigure}
\usepackage{mathabx} 
\numberwithin{equation}{section}
\usepackage{booktabs}
\usepackage{float}
\usepackage{caption}

\newcommand{\rT}{{\rm{T}}}

\newcommand{\dud}[1]{\frac{\partial v}{\partial #1}}
\newcommand{\dudd}[1]{\frac{\partial^2v}{\partial #1^2}}
\newcommand{\duddm}[2]{\frac{\partial^2v}{\partial #1 \partial #2}}
\newcommand{\dd}{\mathrm{d}}

\newcommand{\hV}{{\widehat V}}
\newcommand{\hY}{{\widehat Y}}

\def\C{\mathbb{C}}

\author{Karel~J.~in~'t~Hout\footnote{Department of Mathematics,
		University of Antwerp, Middelheimlaan 1, B-2020 Antwerp, Belgium.
		\mbox{Email}: \texttt{karel.inthout@uantwerpen.be}.}\\\\
		\phantom{\it On the occasion of the 70th birthday of Prof. Peter A. Forsyth}
}

\title{An efficient numerical method for American options and 
their Greeks under the two-asset Kou jump-diffusion model}

\begin{document}
	
\maketitle
	
\begin{abstract}
In this paper we consider the numerical solution of the two-dimensional time-dependent partial 
integro-differential complementarity problem (PIDCP) that holds for the value of American-style 
options under the two-asset Kou jump-diffusion model.	
Following the method of lines (MOL), we derive an efficient numerical method for the pertinent 
PIDCP.
Here, for the discretization of the nonlocal double integral term, an extension is employed of 
the fast algorithm by Toivanen~\cite{T08} in the case of the one-asset Kou jump-diffusion model.
For the temporal discretization, we study a useful family of second-order diagonally implicit 
Runge--Kutta (DIRK) methods.
Their adaptation to the semidiscrete two-dimensional Kou PIDCP is obtained by means of an
effective iteration introduced by d'Halluin, Forsyth \& Labahn~\cite{HFL04} and d'Halluin, 
Forsyth \& Vetzal~\cite{HFV05}.

Ample numerical experiments are presented showing that the proposed numerical method achieves a 
favourable, second-order convergence behaviour to the American two-asset option value as well 
as to its Greeks Delta and Gamma.
\end{abstract}

\section{Introduction}
American-style options are widely traded derivatives in the financial markets due to their
inherent flexibility of early exercise.
For the evolution of the underlying asset prices, jump-diffusion processes constitute a 
principal class of models, see~e.g.~Cont \& Tankov~\cite{CT04book} and Schoutens~\cite{S03book}.
Hence, it is of main importance to know the fair values of American-style options, 
together with their Greeks, under jump-diffusion models.

In this paper we shall deal with American-style options on two underlying assets.
Here a direct generalization of the popular Kou jump-diffusion model for a single asset~\cite{K02} 
to two assets is considered, where the finite-activity jumps in the two asset prices are assumed 
to occur simultaneously and the relative jump sizes possess log-double-exponential distributions.

In general, expressions in (semi-)closed analytical form are lacking in the literature for
the fair values of American options and their Greeks.
Accordingly, there is a big demand for effective approximation methods for these quantities.
In the present paper, we shall consider their approximation via the numerical solution of the 
pertinent two-dimensional time-dependent partial integro-differential complementarity problem 
(PIDCP). 
Here the integral part stems from the contribution of the jumps and is nonlocal: the domain of
integration is the whole first quadrant in the real plane.

For the numerical solution we follow the well-known and versatile method of lines (MOL).
Here the two-dimensional time-dependent PIDCP is first semidiscretized in space (the asset price 
domain) and the resulting large semidiscrete PIDCP system is subsequently discretized in time.
Key elements in the construction of an effective numerical method include the treatment of the 
double integral term, the selection of the temporal discretization scheme, and the handling 
of the early exercise constraint.

In the past two decades much interest has been devoted in the computational finance literature
to the valuation of American options under jump-diffusion models via the numerical solution of 
PIDCPs. 
We give a brief overview of main references in this field, relevant to finite-activity jumps.

d'Halluin, Forsyth \& Labahn~\cite{HFL04} and d'Halluin, Forsyth \& Vetzal~\cite{HFV05} considered
one-dimensional PIDCPs for the values of American options on a single asset under the 
Merton~\cite{M76} and Kou~\cite{K02} jump-diffusion models.
The single integral term has been numerically evaluated by means of the fast Fourier transform, 
cf.~also e.g.~Almendral \& Oosterlee~\cite{AO05},
which avoids the computational burden of directly evaluating matrix-vector products with a 
large dense matrix stemming from the discretization of this term.
For the temporal discretization, the Crank--Nicolson scheme is applied using variable step 
sizes.
Here a fixed-point iteration on the integral part is effectively combined with a penalty 
iteration~\cite{FV02} to handle the early exercise constraint.
Clift \& Forsyth~\cite{CF08} extended the numerical solution approach from \cite{HFL04,HFV05}
to two-dimensional PIDCPs for the values of American options on two assets under generalizations 
of the Merton and Kou jump-diffusion models.

Toivanen~\cite{T08} considered the numerical solution of the one-dimensional Kou PIDCP and 
derived a simple and highly efficient algorithm for evaluating the integral term in this case.
To deal with the early exercise constraint, a useful splitting technique is employed that 
has been introduced and studied by Ikonen \& Toivanen~\cite{IT04,IT09}.

Salmi, Toivanen \& von Sydow~\cite{STS14} proposed, for the temporal discretization of 
two-dimensional PIDCPs, the Crank--Nicolson Adams--Bashforth scheme.
This is a second-order implicit-explicit (IMEX) method that conveniently treats the spatial
differential part in an implicit manner and the integral part in an explicit manner.

Boen \& in~'t~Hout~\cite{BH20} considered the two-dimensional Merton PIDCP and studied, 
for the temporal discretization, a variety of second-order IMEX and alternating direction 
implicit (ADI) methods, using the Ikonen--Toivanen splitting technique to incorporate the 
early exercise constraint.
Here the spatial differential part is again handled implicitly and the integral part explicitly.

Our aim in the present paper is to construct and investigate an efficient numerical method for 
the two-dimensional Kou PIDCP to approximate the fair values of American-style two-asset options 
together with their Greeks Delta and Gamma under the two-asset Kou model.
An outline of the rest of our paper is as follows.

In Section~\ref{Sec_Model}, we formulate the pertinent two-dimensional Kou PIDCP. 

In Section~\ref{Sec_Space}, we consider its semidiscretization.
Here the spatial differential part is discretized in a common way, by applying second-order 
central finite differences on a smooth, nonuniform Cartesian grid.
For the discretization of the double integral term, we employ the two-dimensional extension
recently derived by in~'t~Hout \& Lamotte \cite{HL23} of the fast algorithm by 
Toivanen~\cite{T08} in the case of the one-asset Kou model.
The number of elementary arithmetic operations of this algorithm is directly proportional 
to the number of spatial grid points, which is optimal.

For the temporal discretization of the semidiscrete two-dimensional Kou PIDCP, we investigate
an interesting family of second-order diagonally implicit Runge--Kutta (DIRK) methods 
introduced by Cash~\cite{C04}.
These methods, containing a free parameter $\theta$, are formulated in Section~\ref{Sec_Time}.
For suitable values of $\theta$, the method possesses a strong damping property ($L$-stability). 
This property has been shown to be beneficial for the accurate and stable approximation of the 
Greeks Delta and Gamma for American options in the special case of the Black--Scholes model, 
see Ikonen \& Toivanen~\cite{IT07}, Le Floc'h~\cite{L14} and in~'t~Hout~\cite{H24}.
To our knowledge, the family of DIRK methods has not yet been considered in the literature
for American options and their Greeks under more general, jump-diffusion models.
In the present paper, we shall study its adaptation to the semidiscrete two-dimensional Kou 
PIDCP. 
Here the effective combination of the fixed-point and penalty iterations introduced by d'Halluin 
et al.~\cite{HFL04,HFV05} is employed.

In Section~\ref{Sec_Numer}, ample numerical experiments are presented that provide important 
insight into the convergence behaviour of the proposed numerical method for American option 
values as well as their Greeks Delta and Gamma under the two-asset Kou model.

In Section~\ref{Sec_Conc}, conclusions are given.

\section{American option valuation model}\label{Sec_Model}
Let $v(s_1,s_2,t)$ denote the fair value of an American-style option in the pertinent two-asset Kou 
jump-diffusion model if at $t$ time units before the given maturity time $T$ the two underlying asset 
prices are equal to $s_1\ge 0$ and $s_2\ge 0$.
Let $\phi(s_1,s_2)$ denote the payoff of the option and consider the spatial integro-differential 
operator $\mathcal{A}$ given by
\begin{align}\label{operator}
\mathcal{A}v =
&~ \tfrac{1}{2} \sigma_1^2s_1^2\dudd{s_1} + \rho \sigma_1\sigma_2s_1s_2\duddm{s_1}{s_2} 
+ \tfrac{1}{2} \sigma_2^2 s_2^2\dudd{s_2} + (r-\lambda \zeta_1) s_1\dud{s_1} 
+ (r-\lambda\zeta_2) s_2 \dud{s_2} \nonumber \\
& -(r+\lambda)v+\lambda\int_0^{\infty}\int_0^{\infty} f(y_1,y_2) v(s_1y_1, s_2y_2,t) \dd y_1 \dd y_2.
\end{align}
Here $r$ is the risk-free interest rate, $\sigma_i > 0$ ($i=1,2$) is the instantaneous volatility for 
asset~$i$ conditional on the event that no jumps occur, $\rho$ is the correlation coefficient of 
the two underlying standard Brownian motions, $\lambda$ is the jump intensity of the underlying 
Poisson arrival process and $\zeta_i$ ($i=1,2$) is the expected relative jump size for asset $i$.
The function $f$ is the joint probability density function of two independent random variables having 
log-double-exponential distributions~\cite{K02},
\begin{equation}\label{pdf2D}
f(y_1,y_2) = 
\left\{\begin{array}{lll}
q_1q_2\eta_{q_1}\eta_{q_2}y_1^{\eta_{q_1}-1}y_2^{\eta_{q_2}-1} & (0 < y_1, y_2 < 1),\\
p_1q_2\eta_{p_1}\eta_{q_2}y_1^{-\eta_{p_1}-1}y_2^{\eta_{q_2}-1} & (y_1 \geq 1, 0< y_2 < 1),\\
q_1p_2\eta_{q_1}\eta_{p_2}y_1^{\eta_{q_1}-1}y_2^{-\eta_{p_2}-1} & (0 < y_1 < 1, y_2 \geq 1),\\
p_1p_2\eta_{p_1}\eta_{p_2}y_1^{-\eta_{p_1}-1}y_2^{-\eta_{p_2}-1} & (y_1, y_2 \geq 1).
\end{array}\right.
\end{equation}  
The parameters $p_i$, $q_i$, $\eta_{p_i}$, $\eta_{q_i}$ are all positive constants with $p_i + q_i = 1$ 
and $\eta_{p_i} > 1$.
It holds that
\begin{equation*}
\zeta_i = \frac{p_i\eta_{p_i}}{\eta_{p_i}-1}+\frac{q_i\eta_{q_i}}{\eta_{q_i}+1}-1 \quad (i=1,2).
\end{equation*}
From financial option valuation theory it follows that the function $v$\, satisfies the two-dimensional 
Kou PIDCP
\begin{align}\label{PIDCP}
\begin{cases}
v(s_1,s_2,t) \geq \phi(s_1,s_2),\\\\
\displaystyle\dud{t}(s_1,s_2,t) \geq  \mathcal{A} v(s_1,s_2,t),\\\\
\left(v(s_1,s_2,t)-\phi(s_1,s_2)\right)\left(\displaystyle\dud{t}(s_1,s_2,t)-\mathcal{A} v(s_1,s_2,t)\right) = 0,
\end{cases}
\end{align}
valid pointwise for $(s_1,s_2,t)$ whenever $s_1>0$, $s_2>0$, $0<t\leq T$.
Boundary conditions are given by imposing \eqref{PIDCP} for $s_1=0$ and $s_2=0$, respectively.
The initial condition is provided by the payoff,
\begin{equation}\label{IC}
v(s_1,s_2,0) = \phi (s_1,s_2)
\end{equation}
for $s_1\ge 0$, \,$s_2\ge 0$.

Along with the American option value function $v$, we are interested in this paper in the Greeks 
Delta and Gamma. These important risk quantities are defined by
\begin{equation}\label{2DGreeks}
\Delta_1    = \frac{\partial   v}{\partial s_1}\,,~~ 
\Delta_2    = \frac{\partial   v}{\partial s_2}\,,~~ 
\Gamma_{11} = \frac{\partial^2 v}{\partial s_1^2}\,,~~
\Gamma_{12} = \frac{\partial^2 v}{\partial s_1 \partial s_2}\,,~~
\Gamma_{22} = \frac{\partial^2 v}{\partial s_2^2}\,.
\end{equation}
Clearly, all Greeks \eqref{2DGreeks} appear in the function $\mathcal{A}v$ given by \eqref{operator}.

The three conditions in (\ref{PIDCP}) naturally induce a decomposition of the $(s_1,s_2,t)$-domain: the 
continuation region is the set of all points $(s_1,s_2,t)$ where $\partial v/\partial t= \mathcal{A}v$ 
holds and the early exercise region is the set of all points $(s_1,s_2,t)$ where $v=\phi$ holds.
The joint boundary of these two regions is called the early exercise boundary or free boundary.

We shall deal in this paper with the well-known contract of an American put-on-the-average option, 
which has the payoff
\begin{equation}\label{payoff}
	\phi(s_1,s_2) = \textrm{max} \left(0\,,\,K-\frac{s_1+s_2}{2}\right)
\end{equation}
with given strike price $K>0$.

\section{Spatial discretization}\label{Sec_Space}
As the first step towards the numerical solution of the two-dimensional Kou PIDCP \eqref{PIDCP}, 
we semidiscretize in the spatial variables $(s_1,s_2)$.
The semidiscretization is done in the same way as described in~\cite{HL23} for the case of a 
European put-on-the-average option. 

\begin{figure}
	\centering
	\includegraphics[scale=0.7]{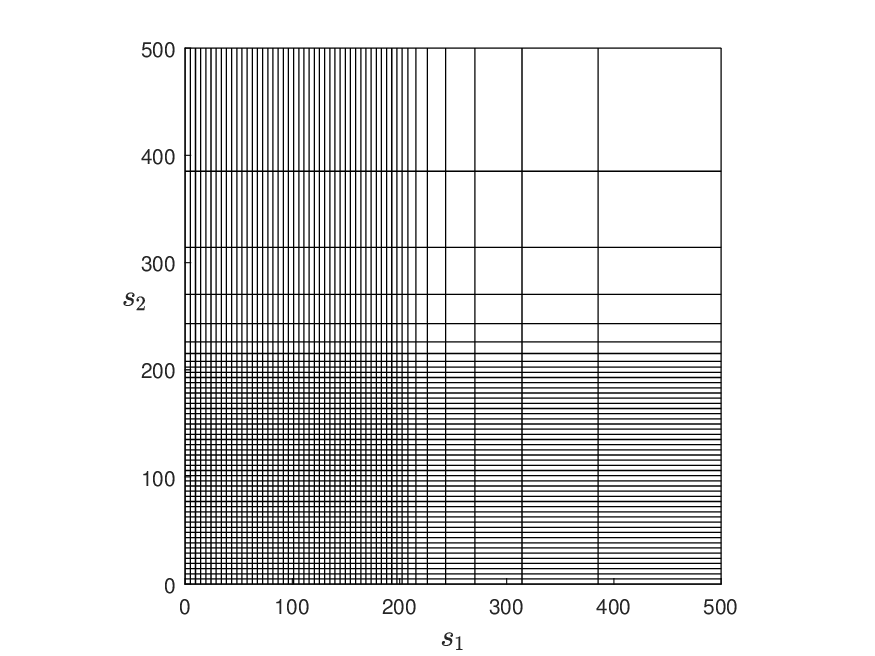}
	\caption{Sample spatial grid for $m=50$, $K=100$, $S_{\rm max} = 500$.}
	\label{grid}
\end{figure}

For any given integer $m\ge 3$, a suitable smooth Cartesian grid $\{ (s_{1,i},s_{2,j}): 0\le i,j \le m \}$ 
is defined in a truncated spatial domain $[0,S_{\rm max}]\times[0,S_{\rm max}]$ with fixed upper bound 
$S_{\rm max} > 2K$ chosen sufficiently large.
In each spatial direction, the mesh is uniform inside the interval $[0,2K]$ and relatively fine.
Outside, it is taken nonuniform and relatively coarse.
Figure~\ref{grid} shows a sample grid where $m=50$, $K=100$ and $S_{\rm max} = 5K$.

Let $V_{i,j}(t)$ denote the semidiscrete approximation to $v(s_{1,i},s_{2,j},t)$ for $0\leq i,j\leq m$ 
and define the vector 
\begin{equation*}
	V(t) = (V_{0,0}(t), V_{1,0}(t),\ldots,V_{m-1,m}(t),V_{m,m}(t))^{\rm T} \in \mathbb{R}^M,
\end{equation*}
where $M=(m+1)^2$ and the symbol $^{\rm T}$ designates the transpose.

All five derivative terms in $\mathcal{A}v$, given by \eqref{operator}, are approximated on the spatial 
grid in a common fashion using second-order central finite differences. 
The obtained semidiscretization of the relevant part of $\mathcal{A}v$, including the $-(r+\lambda)v$ 
term, can be written as $A^{(D)} V(t)$, where $A^{(D)}$ is a given $M\times M$ matrix that is sparse.
Further details are given in~\cite[Sect.~3.1]{HL23}.

For the spatial discretization of the double integral term 
\begin{equation*}\label{integral}
	\mathcal{J} = \lambda\int_0^{\infty}\int_0^{\infty} f(y_1,y_2)v(s_1y_1, s_2y_2,t) \dd y_1 \dd y_2,
\end{equation*}
we employ the two-dimensional extension~\cite{HL23} of the algorithm derived by Toivanen~\cite{T08} 
for the case of the one-asset Kou model.
A simple change of variables yields
\begin{equation*}
	\mathcal{J} = \lambda\int_0^{\infty}\int_0^{\infty} 
	f\bigg(\frac{z_1}{s_1},\frac{z_2}{s_2}\bigg) v(z_1,z_2,t) \frac{\dd z_1 \dd z_2}{s_1 s_2}
\end{equation*}
whenever $s_1, s_2>0$.
Since the density function $f$, given by (\ref{pdf2D}), is defined on a partition of four sets of the 
first quadrant in the real plane, the integral $\mathcal{J}$ is naturally decomposed into four integrals as 
$\mathcal{J} = \mathcal{J}_1 + \mathcal{J}_2 + \mathcal{J}_3 + \mathcal{J}_4$, where
\begin{align*}
	\mathcal{J}_1 &= \lambda q_1q_2\eta_{q_1}\eta_{q_2}s_1^{-\eta_{q_1}}s_2^{-\eta_{q_2}} 
	\int_0^{s_2} \int_0^{s_1} z_1^{\eta_{q_1}-1}z_2^{\eta_{q_2}-1} v(z_1,z_2,t) \dd z_1 \dd z_2,\\
	\mathcal{J}_2 &= \lambda p_1q_2\eta_{p_1}\eta_{q_2}s_1^{\eta_{p_1}}s_2^{-\eta_{q_2}}  
	\int_0^{s_2} \int_{s_1}^{\infty} z_1^{-\eta_{p_1}-1}z_2^{\eta_{q_2}-1} v(z_1,z_2,t) \dd z_1 \dd z_2,\\
	\mathcal{J}_3 &= \lambda q_1p_2\eta_{q_1}\eta_{p_2}s_1^{-\eta_{q_1}}s_2^{\eta_{p_2}}  
	\int_{s_2}^{\infty} \int_0^{s_1} z_1^{\eta_{q_1}-1}z_2^{-\eta_{p_2}-1} v(z_1,z_2,t) \dd z_1 \dd z_2,\\		
	\mathcal{J}_4 &= \lambda p_1p_2\eta_{p_1}\eta_{p_2}s_1^{\eta_{p_1}}s_2^{\eta_{p_2}}   
	\int_{s_2}^{\infty} \int_{s_1}^{\infty} z_1^{-\eta_{p_1}-1}z_2^{-\eta_{p_2}-1} v(z_1,z_2,t) \dd z_1 \dd z_2.
\end{align*}
We shall consider the first integral above. 
For $1\leq i,j \leq m$ let $\mathcal{J}_{1,ij}$ denote its value at the spatial grid point $(s_{1,i},s_{2,j})$.
Write
\begin{equation*}
	\psi_1(s_1,s_2) = \lambda q_1q_2\eta_{q_1}\eta_{q_2}s_1^{-\eta_{q_1}}s_2^{-\eta_{q_2}} 
	\quad \textrm{and} \quad
	\varphi_1(z_1,z_2) = z_1^{\eta_{q_1}-1}z_2^{\eta_{q_2}-1}
\end{equation*}
and for $1\leq k,l \leq m$ define
\begin{equation*}
	\mathcal{G}_{1,kl} = \int_{s_{2,l-1}}^{s_{2,l}} \int_{s_{1,k-1}}^{s_{1,k}} 
	\varphi_1(z_1,z_2) v(z_1,z_2,t) \dd z_1 \dd z_2.
\end{equation*}
Then $\mathcal{J}_{1,ij}$ can be expressed via a double cumulative sum,
\begin{equation}\label{cumsum1}
	\mathcal{J}_{1,ij} = \psi_1(s_{1,i},s_{2,j}) \sum_{k=1}^i \sum_{l=1}^j \mathcal{G}_{1,kl}
	\quad
	(1\leq i, j \leq m).
\end{equation}
The obvious but important fact holds that the $\mathcal{G}_{1,kl}$ are independent of the indices $i$ and $j$.
Consequently, if all values $\mathcal{G}_{1,kl}$ are given, then computing the double cumulative sums in 
\eqref{cumsum1} for all $i$, $j$ can be performed, to leading order, in just $2m^2$ additions.

For any given $k, l$ with $1\leq k, l \leq m$ an approximation $G_{1,kl}$ to $\mathcal{G}_{1,kl}$ is obtained
by using the bilinear interpolant ${\widetilde v}_{kl}$ that approximates $v$ on the $(z_1,z_2)$-domain 
$[s_{1,k-1}, s_{1,k}] \times [s_{2,l-1}, s_{2,l}]$:
\begin{equation*}
	G_{1,kl} = \int_{s_{2,l-1}}^{s_{2,l}} \int_{s_{1,k-1}}^{s_{1,k}} 
	\varphi_1(z_1,z_2) {\widetilde v}_{kl}(z_1, z_2,t) \dd z_1 \dd z_2.
\end{equation*}
The approximation $J_{1,ij}$ to $\mathcal{J}_{1,ij}$ is then defined by
\begin{equation}\label{cumsum1_semi}
	J_{1,ij} = \psi_1(s_{1,i},s_{2,j}) \sum_{k=1}^i \sum_{l=1}^j G_{1,kl}
	\quad
	(1\leq i, j \leq m).
\end{equation}
It is easily verified that one arrives at the simple formula
\begin{equation}\label{G1kl}
	G_{1,kl} = 
	\gamma_{1,kl}^{00} V_{k-1,l-1}(t) + \gamma_{1,kl}^{10} V_{k,l-1}(t) +
	\gamma_{1,kl}^{01} V_{k-1,l}(t)   + \gamma_{1,kl}^{11} V_{k,l}(t) 
\end{equation}
with certain real coefficients $\gamma_{1,kl}^{ab}$ for $a,b\in \{0,1\}$ that are fully determined by the Kou 
parameters and the spatial grid (whenever $1\leq k, l \leq m$).
These coefficients are independent of $t$ and they can therefore be computed upfront, before the time discretization.
Since the values of $\psi_1$ in \eqref{cumsum1_semi} can also be computed upfront, it readily follows that 
the number of elementary arithmetic operations to compute, for any given $t$, all approximations $J_{1,ij}$ 
($1\leq i,j\leq m$) by \eqref{cumsum1_semi} is, to leading order, equal to $10m^2$.

For the other three integrals $\mathcal{J}_2, \mathcal{J}_3, \mathcal{J}_4$, the discretization and efficient 
evaluation is performed completely analogous to that for $\mathcal{J}_1$.
The favourable result is obtained that the number of elementary arithmetic operations to compute, for any given $t$, 
the approximation to $\mathcal{J}$ on the full spatial grid is, to leading order, directly proportional to the 
number of spatial grid points, which is optimal.
For further details we refer to~\cite[Sect.~3.2]{HL23}.

For notational convenience, the approximation to the double integral $\mathcal{J}$ on the full spatial grid will 
formally be denoted by $A^{(J)} V(t)$, where $A^{(J)}$ is a given $M\times M$ matrix.
We emphasize that in our paper matrix-vector products in the case of $A^{(J)}$ are never evaluated directly, since 
this matrix is large and dense, but are always computed by the efficient algorithm described above.

Write $A = A^{(D)} + A^{(J)}$. 
Then the spatial discretization gives rise to a large semidiscrete PIDCP system of the form
\begin{equation}\label{lcp_ode}
V(t) \ge V^0,
\quad V'(t) \ge AV(t),
\quad (V(t) - V^0)^\rT (V'(t) - AV(t))=0
\end{equation}
for $0 < t \le T$ with $V(0)=V^0$.
Here vector inequalities are to be understood componentwise.
The initial vector $V^0\in \mathbb{R}^M$ is given by pointwise evaluation of the payoff function $\phi$ on 
the spatial grid, except that cell averaging is used near the line segment $s_1+s_2 = 2K$, where $\phi$ is 
nonsmooth, see~e.g.~\cite{H17book,HL23}.

Semidiscrete approximations to all Greeks Delta and Gamma are directly acquired by application of the 
pertinent second-order central finite difference formulas using the entries from the vector $V(t)$.
Hence, as such, they form terms in the matrix-vector product $A^{(D)} V(t)$.

\section{Temporal discretization}\label{Sec_Time}
For the temporal discretization of the semidiscrete PIDCP system \eqref{lcp_ode}, we shall investigate in 
this paper the adaptation of the {\it diagonally implicit Runge--Kutta (DIRK) method} given by the Butcher 
tableau
\begin{equation}\label{tableau_DIRK}
\renewcommand\arraystretch{1.2}
\begin{array}
{c|ccc}
0 & 0 \\
1 & 1-\theta & \theta\\
1 &\frac{1}{2} & \frac{1}{2}-\theta & \theta\\
\hline
&\frac{1}{2} & \frac{1}{2}-\theta & \theta
\end{array}
\end{equation}
where $\theta >0$ denotes a given parameter.
Method \eqref{tableau_DIRK} possesses classical order of consistency\footnote{That is, for fixed systems 
of ordinary differential equations.} equal to two for any value of $\theta$ and has the stability function 
given by
\begin{equation}\label{stabfunc_DIRK}
R(z) = \frac{1+(1-2\theta)z+(\tfrac{1}{2}-2\theta+\theta^2)z^2}{(1-\theta z)^2} \quad (z\in \C).
\end{equation}
It is well-known that \eqref{tableau_DIRK} is $A$-stable\footnote{A Runge--Kutta method is called 
{\it $A$-stable} if its stability function $R$ satisfies $|R(z)|\le 1$ whenever $z\in \C$, $\Re z \le 0$. 
If in addition $R(\infty)=0$, then the method is called {\it $L$-stable}, see e.g.~\cite{HW91,HV03}.} 
for all $\theta\ge \tfrac{1}{4}$ and $L$-stable if and only if $\theta = 1 \pm \tfrac{1}{2} \sqrt{2}$.

The DIRK method \eqref{tableau_DIRK} was introduced by Cash~\cite{C04} and has been considered by 
various authors in the computational finance literature, notably for American option valuation.
Khaliq, Voss \& Kazmi~\cite{KVK06} studied an adaptation of this method for the numerical valuation of 
American one- and two-asset options under the Black--Scholes model. 
Ikonen \& Toivanen~\cite{IT07,IT09} examined adaptations for the numerical valuation of American options 
under the one-asset Black--Scholes and Heston models.
Recently, in~'t~Hout~\cite{H24} investigated an adaptation of \eqref{tableau_DIRK} for the approximation
of the Greeks Delta and Gamma in the case of American one- and two-asset options under the Black--Scholes 
model.

A common choice for the parameter value in \eqref{tableau_DIRK} is \mbox{$\theta = 1-\tfrac{1}{2}\sqrt{2}$}, 
which yields both $L$-stability and a relatively small error constant.
This choice has been considered in all references above.
It is worth mentioning that for this value of $\theta$ the stability function \eqref{stabfunc_DIRK} is 
identical to that of the TR-BDF2 method, a familiar combination of the trapezoidal rule and the 
second-order backward differentiation formula introduced by Bank et al.~\cite{BCF85} and subsequently 
studied by e.g.~Hosea \& Shampine~\cite{HS96}.
The TR-BDF2 method has been advocated for the numerical approximation of one-asset American option 
values and their Greeks Delta and Gamma by Le Floc'h~\cite{L14}.

It is further interesting to note that \eqref{tableau_DIRK} forms the underlying implicit method of both 
the Hundsdorfer--Verwer (HV) scheme and the modified Craig--Sneyd (MCS) scheme, two well-known alternating 
direction implicit (ADI) schemes that have first been studied for the temporal discretization of (semidiscrete)
PDEs in finance by in~'t~Hout and Welfert~\cite{HW07,HW09}.
For the MCS scheme, the value $\theta = \tfrac{1}{3}$ is often judiciously selected, cf.~e.g.~\cite{H24,HL23}.
Then the DIRK method \eqref{tableau_DIRK} is $A$-stable, but not $L$-stable, since $R(\infty)=-\tfrac{1}{2}$.

For the adaptation of \eqref{tableau_DIRK} to \eqref{lcp_ode}, we shall adopt the popular penalty iteration. 
This has first been studied in the literature on American option valuation by Zvan, Forsyth \& Vetzal~\cite{ZFV98,ZFV01} 
and Forsyth \& Vetzal~\cite{FV02} and has since then been widely considered in the computational finance literature.
Along with it, a fixed-point iteration is employed on the integral part.
This useful technique has first been proposed and analyzed in the literature on European option valuation under 
jump-diffusion models by Tavella \& Randall~\cite{TR00book} and d'Halluin, Forsyth \& Vetzal~\cite{HFV05}, 
respectively.
A simple and effective combination of the penalty iteration and the fixed-point iteration has been studied by 
d'Halluin, Forsyth \& Labahn~\cite{HFL04} and Clift \& Forsyth~\cite{CF08} for American option valuation 
under jump-diffusion models.
Following their approach, we acquire the adaptation of the DIRK method \eqref{tableau_DIRK} under consideration 
in the present paper.

Let ${\it Large} >0$ be a given, fixed large number and let ${\it tol} >0$ be a given, fixed small tolerance.
Let $0=t^0<t^1<t^2<\ldots<t^N=T$ be any given sequence of temporal grid points with (constant or variable) 
step sizes $\Delta t^n = t^n-t^{n-1}$ ($1\le n\le N$) and set $\hV^{0} = V^0$.
Then the adaptation of method \eqref{tableau_DIRK} to the semidiscrete PIDCP \eqref{lcp_ode} generates 
approximations $\hV^{n}$ to $V(t^n)$ successively, in a one-step manner, for $n=1,2,3,\ldots,N$ by the 
{\it DIRK-P method}\,:
\begin{equation}\label{DIRK_P}
\left\{\begin{array}{l}
W_1 = \hV^{n-1} +(1-\theta)\Delta t^n (A^{(D)}\hV^{n-1} + A^{(J)} \hV^{n-1}), \\\\
\left( I-\theta\Delta t^n A^{(D)} +P_{k-1} \right) Y_k =  W_1 + \theta\Delta t^n A^{(J)} Y_{k-1} + P_{k-1} V^0 \quad (k=1,2,\ldots,\kappa_1), \\\\
\hY = Y_{\kappa_1}\,, \\\\
W_2 = \hV^{n-1}+\tfrac{1}{2}\Delta t^n (A^{(D)}\hV^{n-1} + A^{(J)}\hV^{n-1}) + 
\left(\tfrac{1}{2}-\theta\right)\Delta t^n (A^{(D)}\hY + A^{(J)}\hY), \\\\
\left( I-\theta\Delta t^n A^{(D)} +Q_{k-1} \right) Z_k = W_2 + \theta\Delta t^n A^{(J)} Z_{k-1} + Q_{k-1} V^0  \quad (k=1,2,\ldots,\kappa_2), \\\\
\hV^{n} = Z_{\kappa_2}\,.
\end{array}\right.
\end{equation}
In each time step, two consecutive iteration processes are performed. 
Here $P_k$, respectively $Q_k$, is defined as the diagonal matrix with $l$-th diagonal entry equal to 
${\it Large}$ if $Y_{k,l} < V^0_l$, respectively $Z_{k,l} < V^0_l$, and zero otherwise ($k\ge 0$).
For the first iteration process we consider the natural stopping criterion
\begin{equation*}\label{conv_crit}
\max_l \frac{| Y_{\kappa_1,l}-Y_{\kappa_1-1,l}|}{\max \{1,| Y_{\kappa_1,l} | \} } < {\it tol} 
\quad {\rm or} \quad P_{\kappa_1} = P_{\kappa_1-1},
\end{equation*}
and similarly for the second iteration process.
For the solution of the large, sparse linear systems in \eqref{DIRK_P}, the BiCGSTAB iterative 
method\footnote{As implemented in Matlab version R2020b through the function {\tt bicgstab}.} is used 
with an ILU preconditioner.

Throughout this paper the common values ${\it Large} = 10^{7}$ and ${\it tol} = 10^{-7}$ are selected.
As starting vectors for the two iteration processes, the standard choice in the literature is
$Y_0 = Z_0 = \hV^{n-1}$. 
In the present paper, we shall consider linear extrapolation from the temporal grid points 
$t^{n-2}$ and $t^{n-1}$ to $t=t^n$ to define these vectors (whenever $n\ge 2$)\,:
\begin{equation*}
Y_0 = Z_0 = \hV^{n-1} + \frac{\Delta t^n}{\Delta t^{n-1}} (\hV^{n-1}-\hV^{n-2}).
\end{equation*}
In our numerical experiments this choice is found to reduce the total number of iterations compared 
to the standard choice (which can be viewed as constant extrapolation).

Analogously to Section~\ref{Sec_Space}, fully discrete approximations to all Greeks Delta and Gamma 
are directly acquired by application of the pertinent second-order central finite difference formulas 
using the entries from the vector $\hV^{n}$.
Hence, they are constituents of the matrix-vector product $A^{(D)} \hV^{n}$.

We remark that the first two lines of \eqref{DIRK_P} define the temporal discretization method 
that has been studied in \cite{CF08,HFL04} and alluded to above (upon setting $\hV^{n} = \hY$).
In this case, the underlying Runge--Kutta method is the classical $\theta$-method.

It is well-known in the literature on American option valuation that suitable variable step sizes can 
lead to an improved convergence behaviour of temporal discretization methods compared to the use of 
constant step sizes, see e.g.~Forsyth \& Vetzal~\cite{FV02}, Ikonen \& Toivanen~\cite{IT09} and 
Reisinger \& Whitley~\cite{RW14}.
In this paper, we shall consider the nonuniform temporal grid defined by~\cite{IT09,RW14}
\begin{equation}\label{variablestep}
t^n = \left( \frac{n}{N} \right)^2 T \quad {\rm for}~ n=0,1,2,\ldots,N.
\end{equation}
Here, the variable step size $\Delta t^n$ is smallest for $n=1$ and grows linearly with~$n$.

\section{Numerical study}\label{Sec_Numer}
In this section we consider four instances of the DIRK-P method \eqref{DIRK_P}:\\\\
\begin{tabular}{lll}
\mbox{DIRKa-P}\,:& $\theta=1-\tfrac{1}{2}\sqrt{2}$ \\
\mbox{DIRKb-P}\,:& $\theta=\tfrac{1}{3}$ \\ 
\mbox{DIRKc-P}\,:& $\theta=1$            \\ 
\mbox{DIRKd-P}\,:& $\theta=1+\tfrac{1}{2}\sqrt{2}$  
\end{tabular}
\\\\
For each of these four instances, the underlying Runge--Kutta method is $A$-stable.
For DIRKa-P and DIRKd-P it is also $L$-stable, whereas for DIRKb-P and DIRKc-P it 
is not (then $R(\infty)=-\tfrac{1}{2}$).
In view of this, we apply the latter two methods with backward Euler damping, i.e.,
the first two time steps ($n = 1,2$) are replaced by the BE-P method, which is defined
by the first two lines of \eqref{DIRK_P} with $\theta = 1$ and setting $\hV^{n} = \hY$.

We shall investigate the {\it temporal discretization errors}\footnote{For this concept, see 
e.g.~Hundsdorfer \& Verwer~\cite{HV03}.} of the above four DIRK-P methods at $t=t_N=T$ on a 
given {\it region of interest} ROI in the $(s_1,s_2)$-domain, defined by
\begin{equation}\label{temp_error}
\widehat{E}(m,N) = 
\max \left\{ |V_{i,j}(T)-\hV^N_{i,j}|:\, 0\le i,j\le m,~(s_{1,i}\,,\, s_{2,j}) \in \textrm{ROI}\, \right\}.
\end{equation}
Here $V(T)\in\mathbb{R}^M$ represents\footnote{A reference solution has been computed for each
pertinent $m$ by applying the DIRKa-P method with $N=500$ time steps.} the exact solution to 
the semidiscrete PIDCP system \eqref{lcp_ode} at $t=T$.
Along with \eqref{temp_error}, the temporal discretization errors in the case of the Greeks 
Delta and Gamma will be studied.
These errors are defined completely analogously to \eqref{temp_error} and are denoted by
$\widehat{E}_{\Delta_i} (m,N)$ and $\widehat{E}_{\,\Gamma_{ij}} (m,N)$ whenever $1\le i,j \le 2$.
Notice that all temporal discretization errors are measured in the important maximum-norm.

For the numerical experiments we consider a typical parameter set, given by Table~\ref{paramset}.
It forms a blend of the parameter values for the Black--Scholes model in~\cite{H24} and the Kou 
model in~\cite{CF08}.
We take $S_{\rm max} = 10K$.

\begin{table}[H]
\centering
\caption{Parameter set for the two-asset Kou jump-diffusion model.}
\begin{tabular}{@{}cccccccccccccc@{}}
\toprule
& $\sigma_1$ & $\sigma_2$ & $r$ & $\rho$ & $\lambda$ & $p_1$ & $p_2$ & $\eta_{p_1}$ & $\eta_{q_1}$ & $\eta_{p_2}$ & $\eta_{q_2}$ & $K$ & $T$ \\ 
\midrule
& 0.30 & 0.40 & 0.01 & 0.50 & 0.50 & 0.40 & 0.60 & $1/0.20$ & $1/0.15$ & $1/0.18$ & $1/0.14$ & 100 & 0.5 \\
\bottomrule
\end{tabular}
\label{paramset}
\end{table}

\begin{figure}[H]
	\centering
	\includegraphics[scale=0.49]{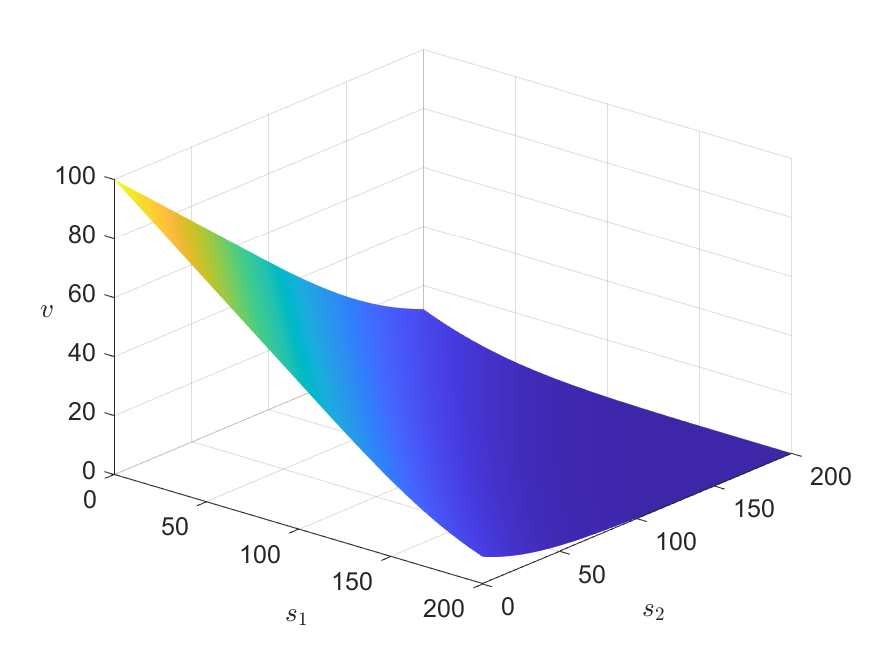}
	\includegraphics[scale=0.49]{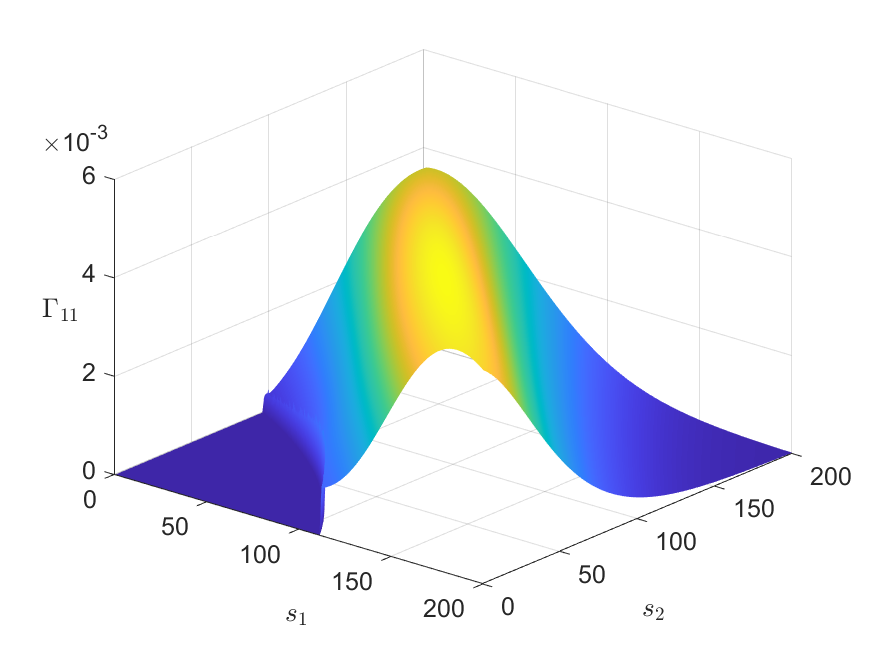}\\
	\includegraphics[scale=0.49]{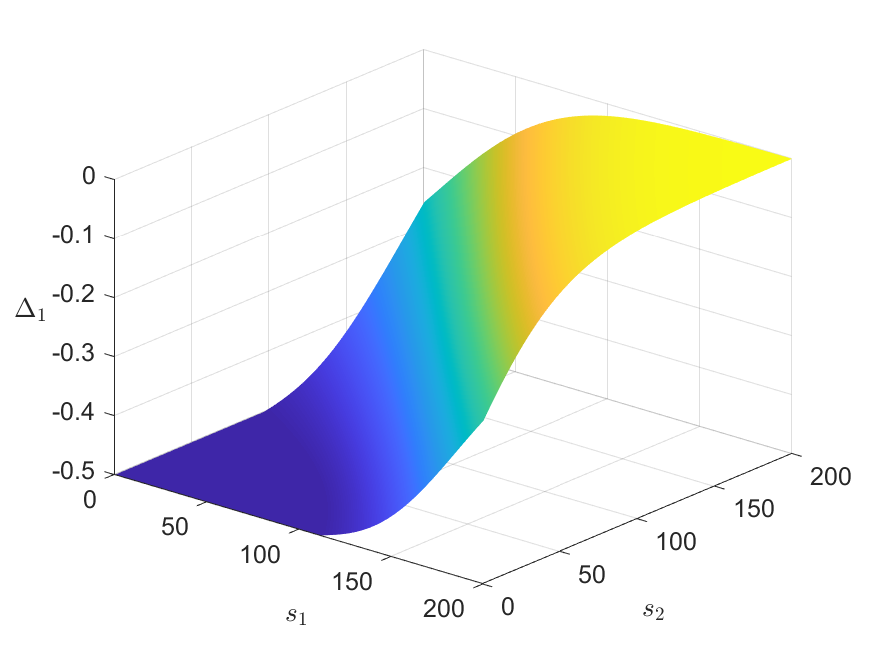}
	\includegraphics[scale=0.49]{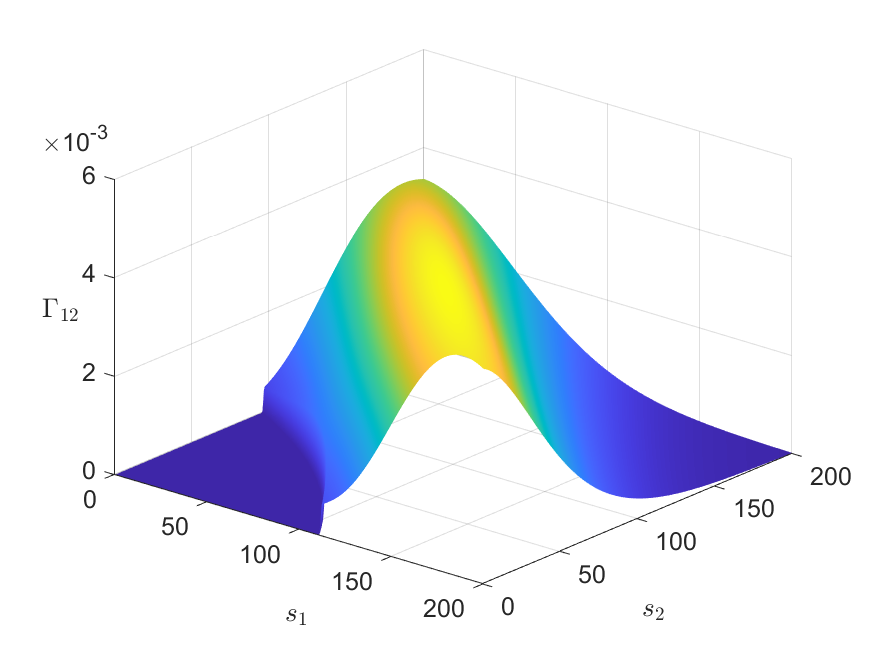}\\
	\includegraphics[scale=0.49]{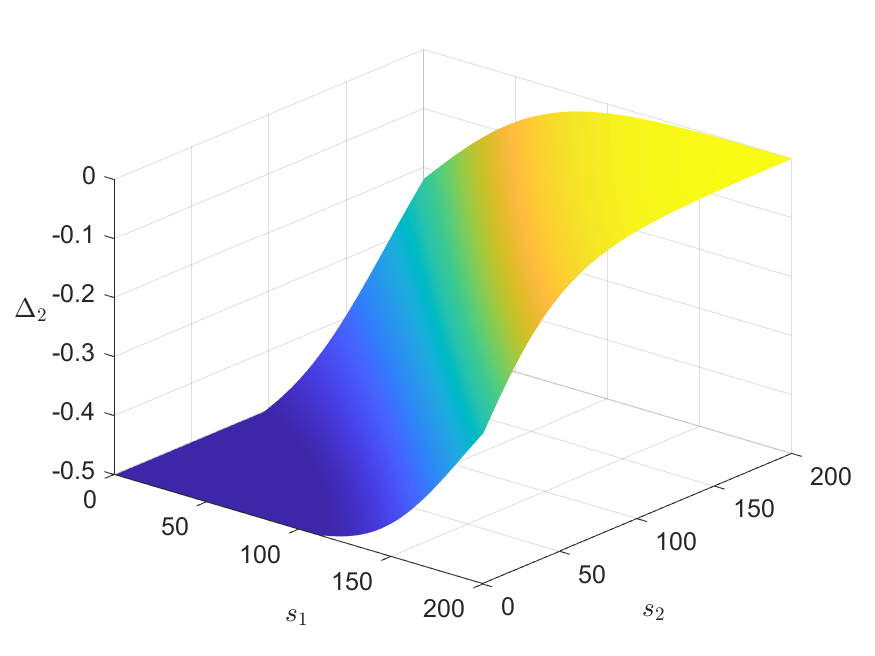}
	\includegraphics[scale=0.49]{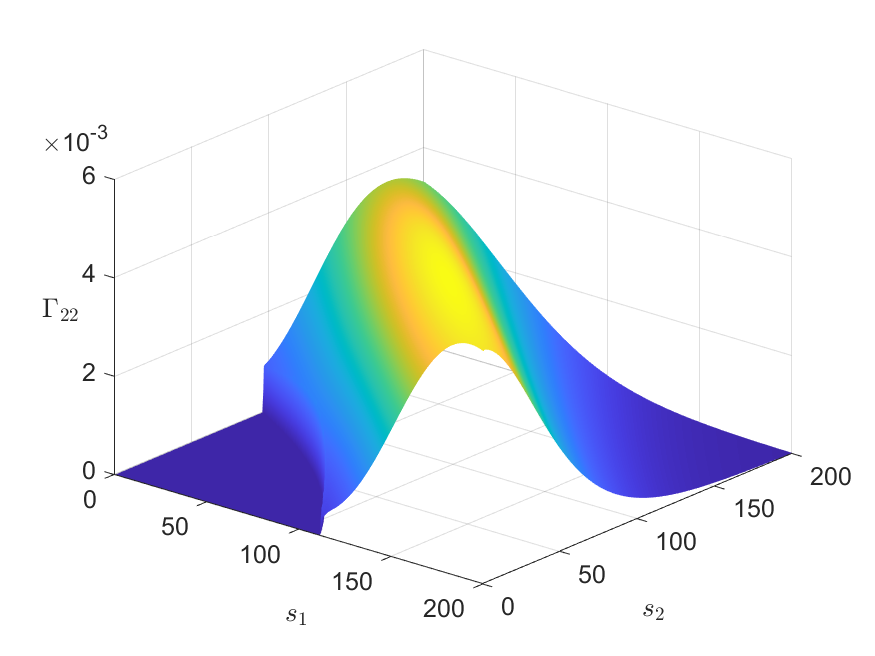}\\
	\caption{The value, Delta and Gamma functions of an American put-on-the-average option for $t=T$ 
	and parameter set given by Table~\ref{paramset}.}
	\label{graphs}
\end{figure}

\begin{figure}
	\centering
	\includegraphics[scale=0.6]{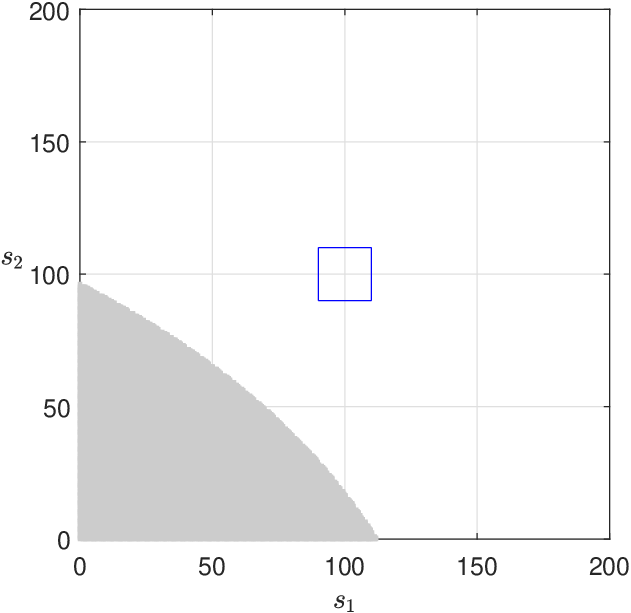}
	\caption{In grey: the early exercise region of an American put-on-the-average option for $t=T$ 
	and parameter set given by Table~\ref{paramset}. In blue: the region of interest.}
	\label{EER}
\end{figure}

Figure~\ref{graphs} shows the approximated graphs of the American put-on-the-average option 
value function and the five Delta and Gamma functions~\eqref{2DGreeks} for $t=T$ on the 
$s$-domain $[0,2K]\times [0,2K]$.
Figure~\ref{EER} shows in grey the corresponding early exercise region for $t=T$.
Here the blue square indicates the region of interest that we consider, given by
$\textrm{ROI} = (0.9K, 1.1K)\times (0.9K, 1.1K)$.
Clearly, this ROI lies well within the continuation region, at a significant distance from the 
early exercise boundary.

\begin{figure}[H]
	\centering
	\includegraphics[scale=0.5]{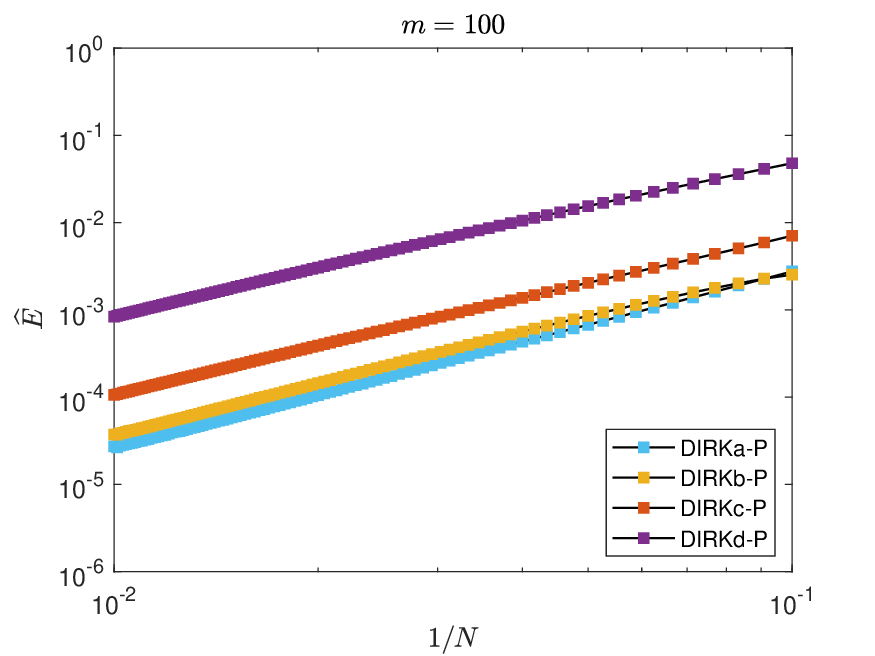}
	\includegraphics[scale=0.5]{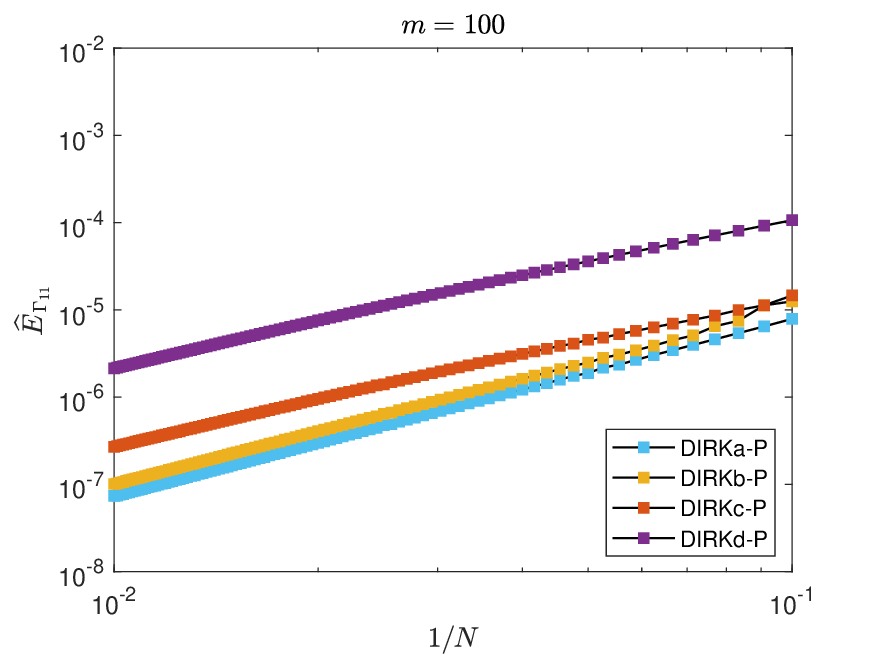}\\
	\includegraphics[scale=0.5]{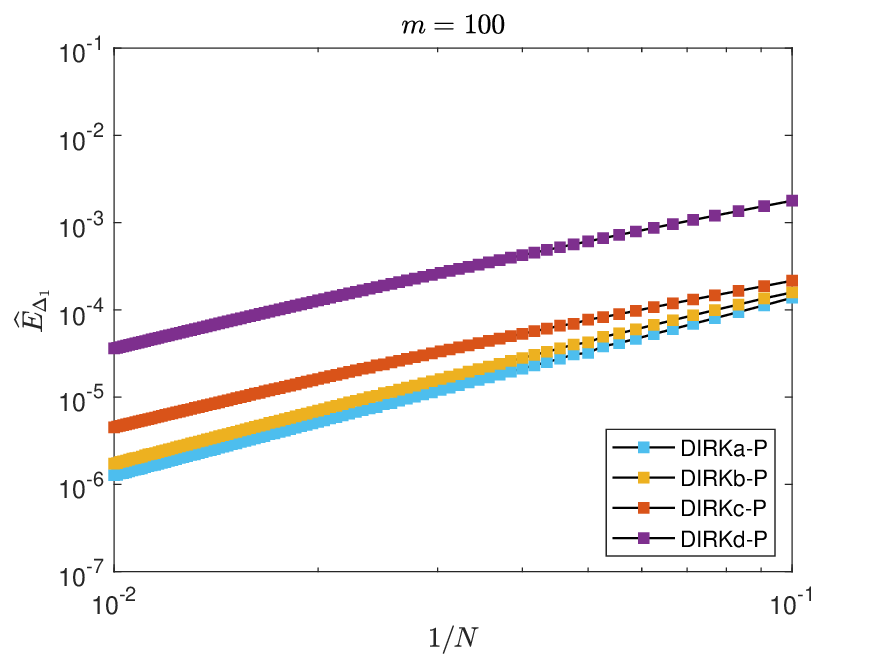}
	\includegraphics[scale=0.5]{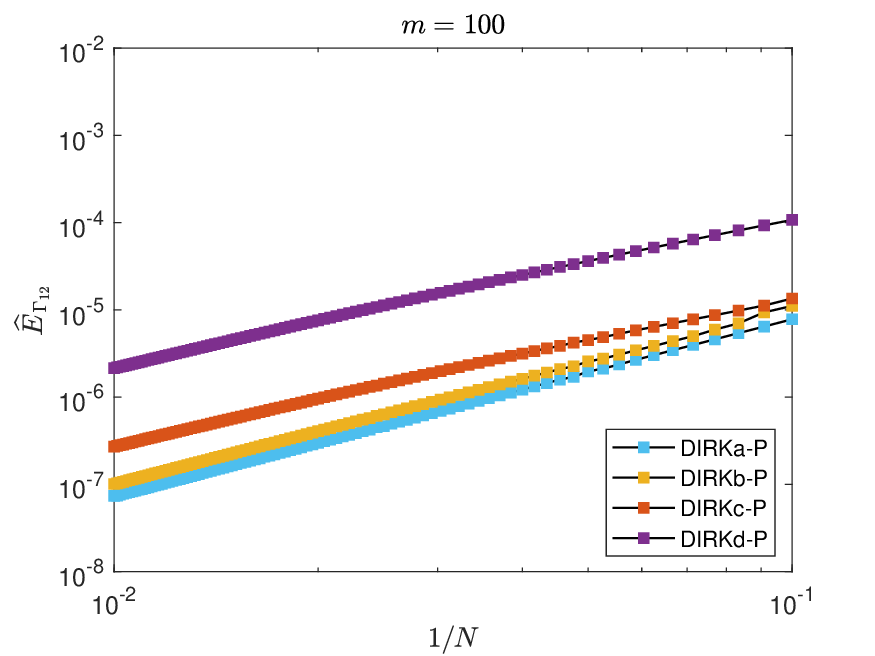}\\
	\includegraphics[scale=0.5]{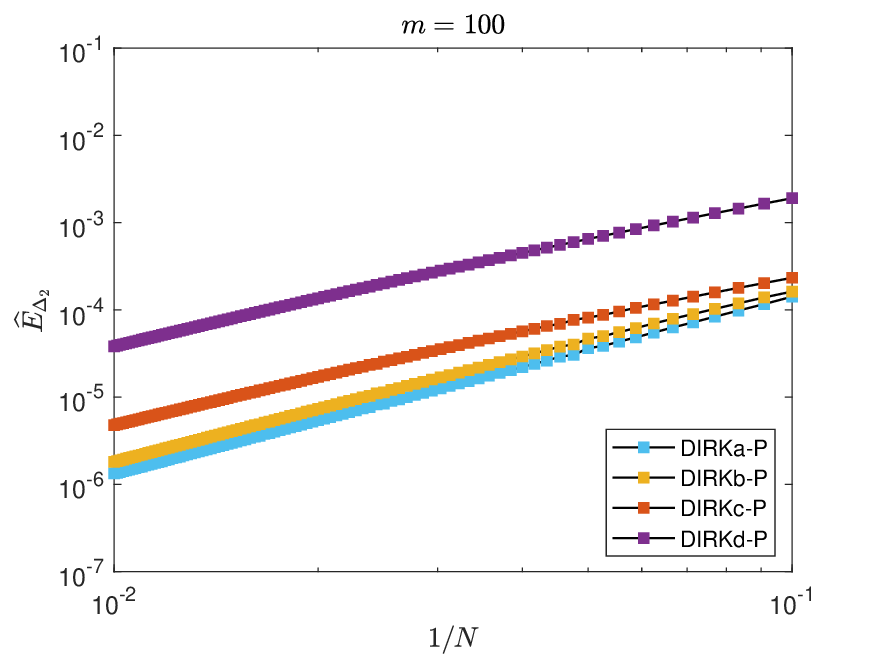}
	\includegraphics[scale=0.5]{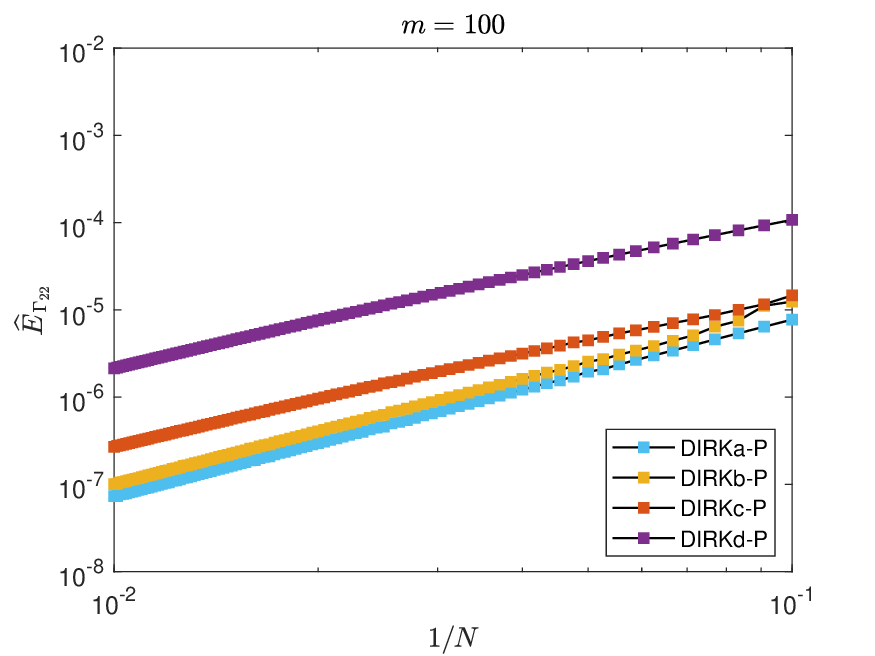}\\
	\caption{ Temporal discretization errors of the DIRKa-P, DIRKb-P, DIRKc-P, DIRKd-P methods 
	for $m=100$.
	Option value (top left), $\Delta_1$ (middle left), $\Delta_2$ (bottom left), $\Gamma_{11}$ 
	(top right), $\Gamma_{12}$ (middle right), $\Gamma_{22}$ (bottom right).}	
	\label{2Derror1}
\end{figure}

\begin{figure}[H]
	\centering
	\includegraphics[scale=0.5]{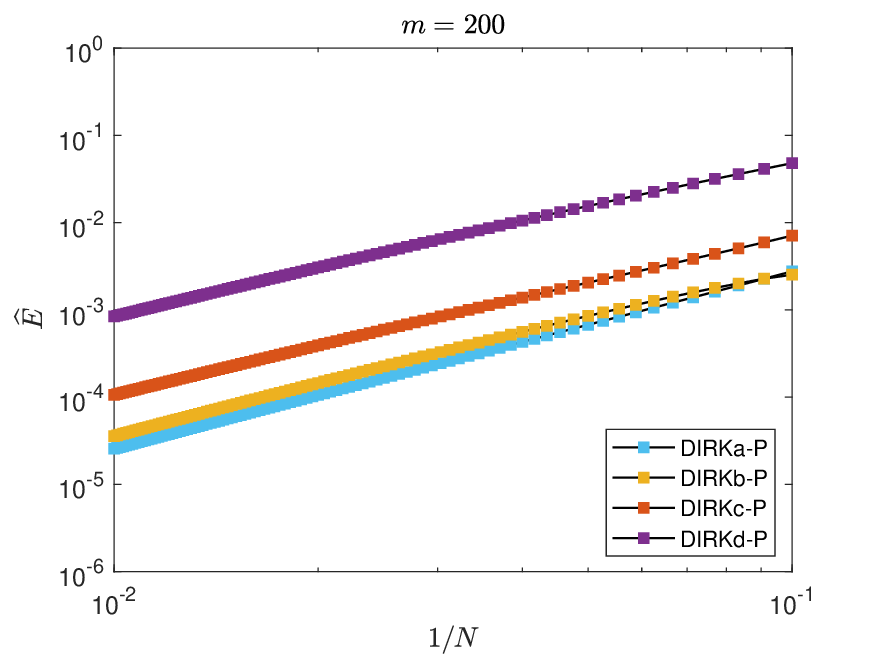}
	\includegraphics[scale=0.5]{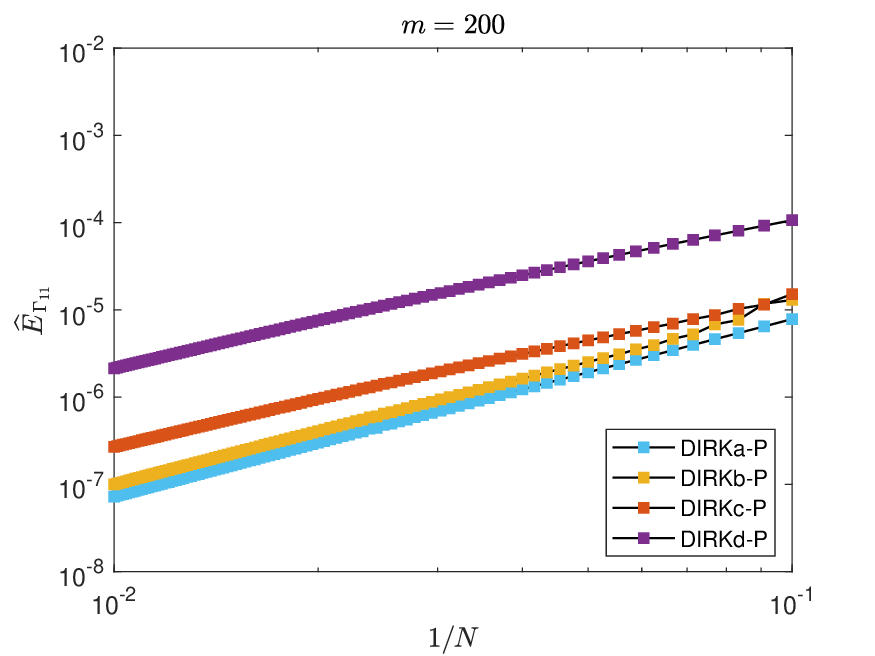}\\
	\includegraphics[scale=0.5]{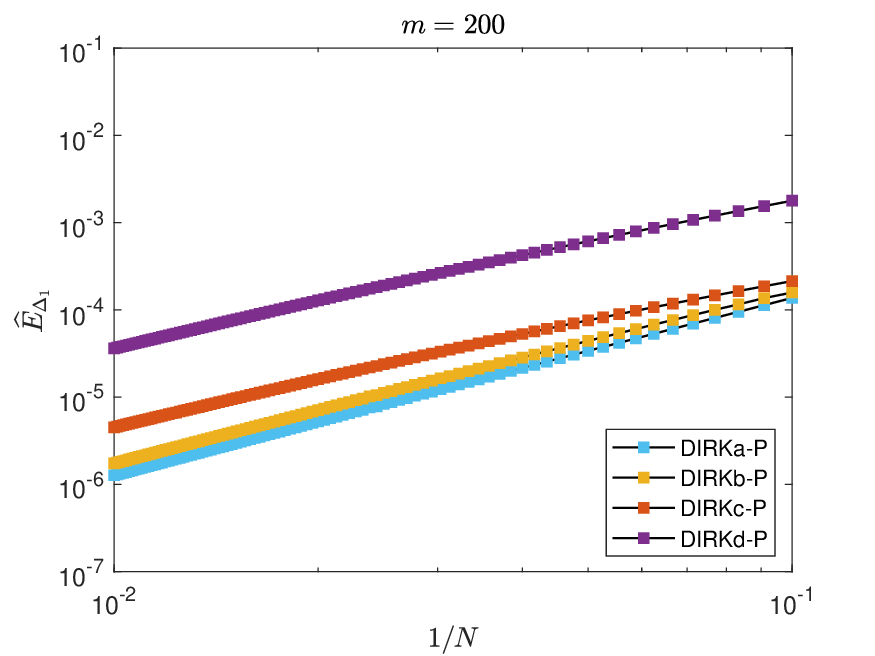}
	\includegraphics[scale=0.5]{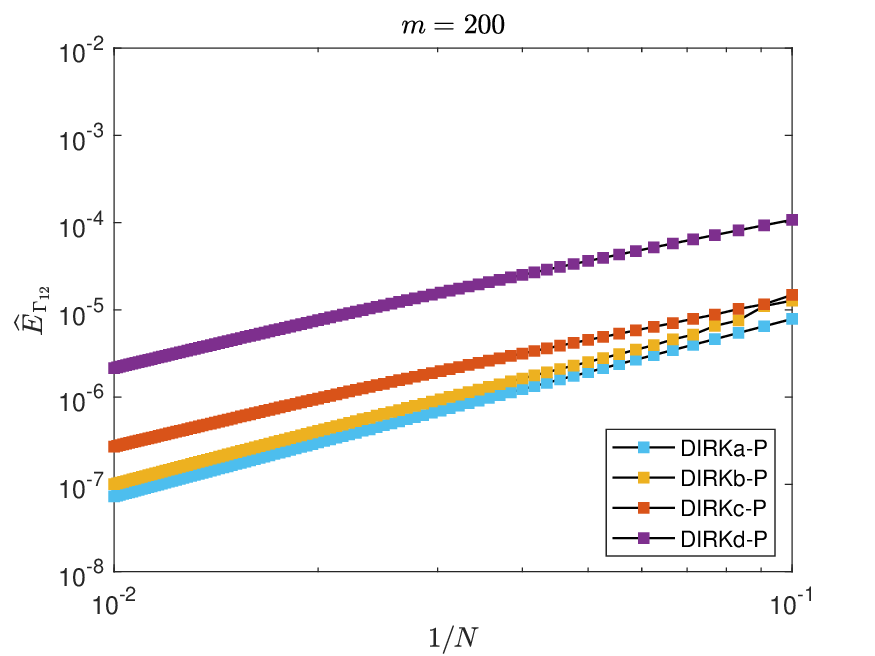}\\
	\includegraphics[scale=0.5]{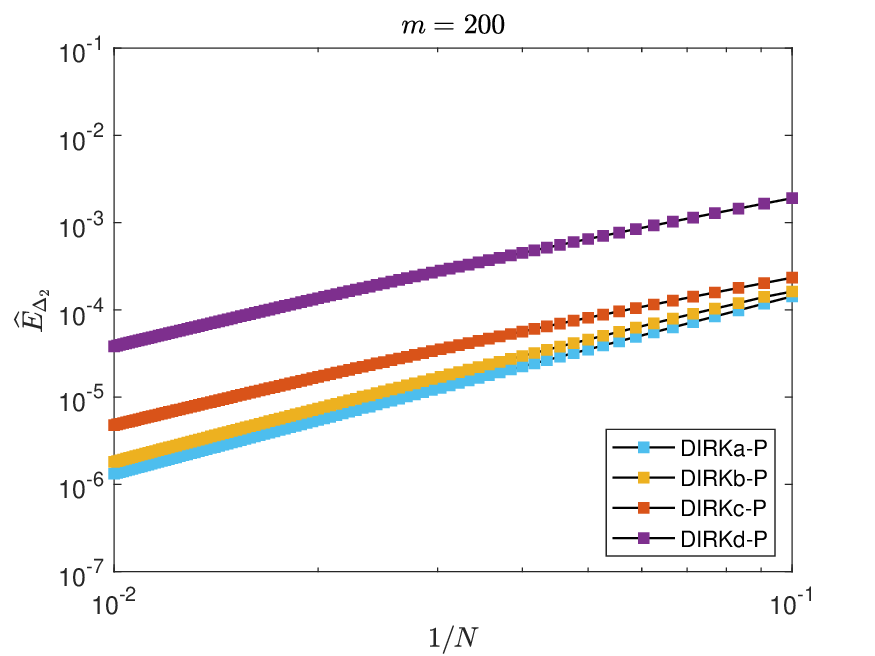}
	\includegraphics[scale=0.5]{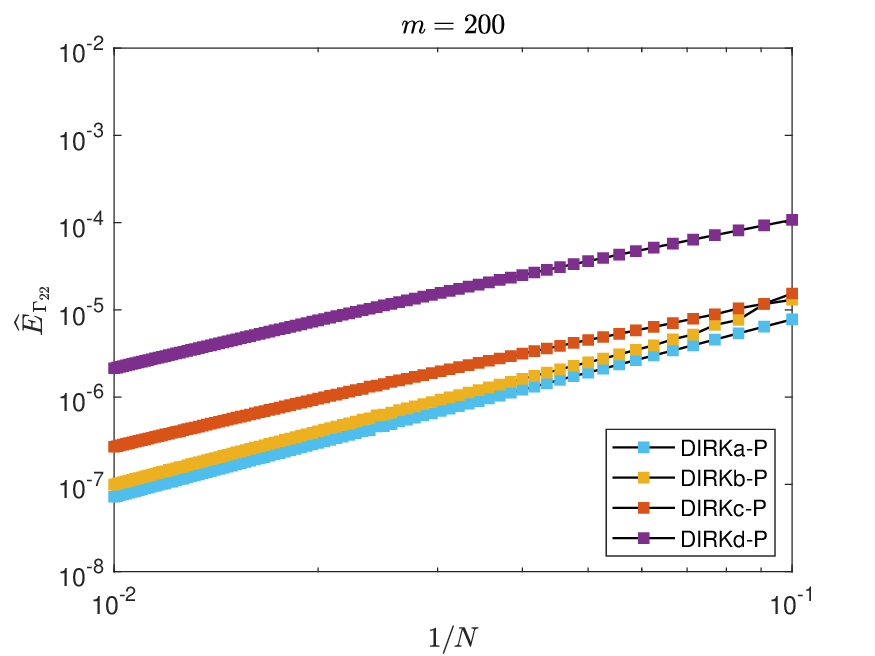}\\
	\caption{ Temporal discretization errors of the DIRKa-P, DIRKb-P, DIRKc-P, DIRKd-P methods 
	for $m=200$.
	Option value (top left), $\Delta_1$ (middle left), $\Delta_2$ (bottom left), $\Gamma_{11}$ 
	(top right), $\Gamma_{12}$ (middle right), $\Gamma_{22}$ (bottom right).}	
	\label{2Derror2}
\end{figure}

\begin{figure}[H]
	\centering
	\includegraphics[scale=0.5]{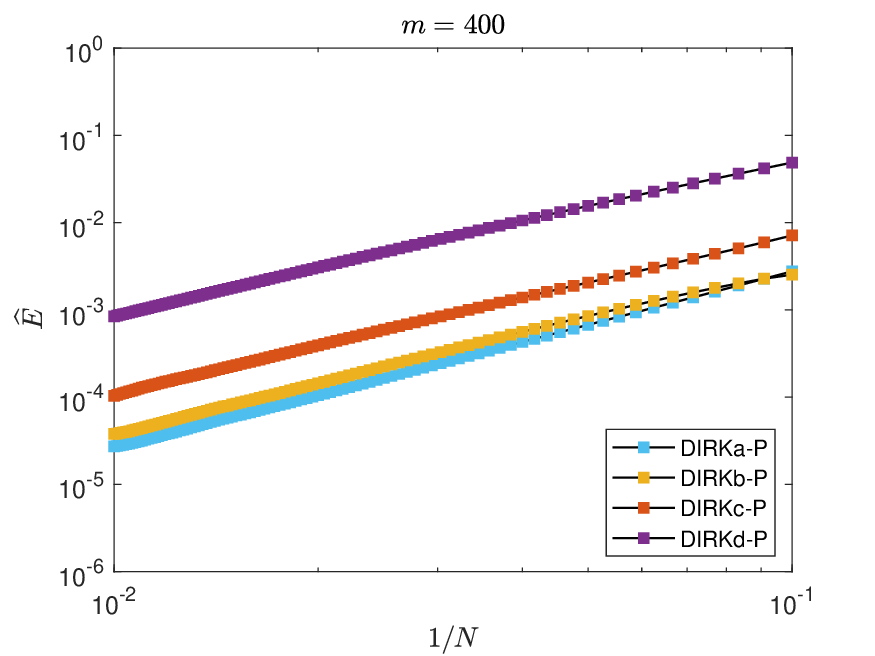}
	\includegraphics[scale=0.5]{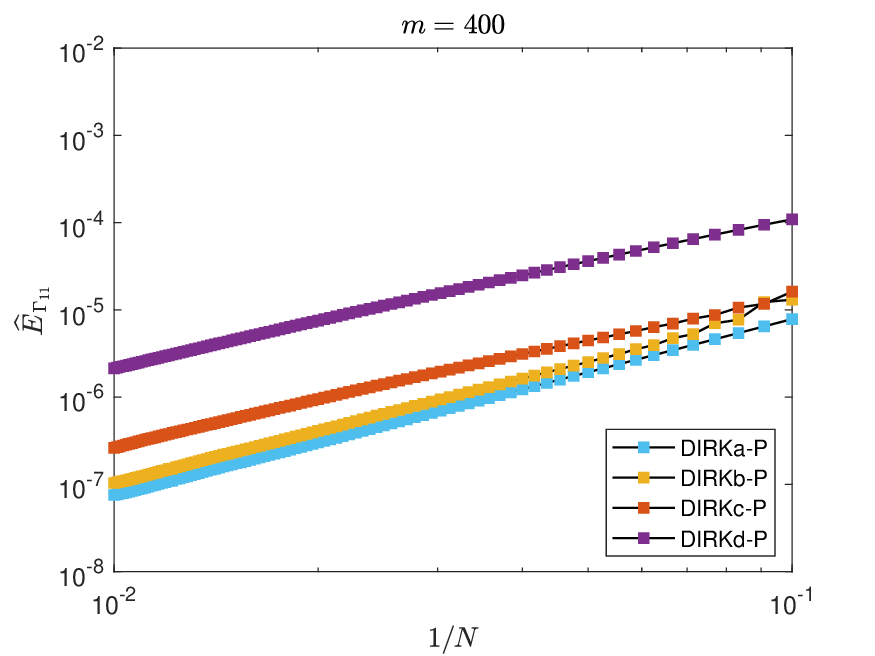}\\
	\includegraphics[scale=0.5]{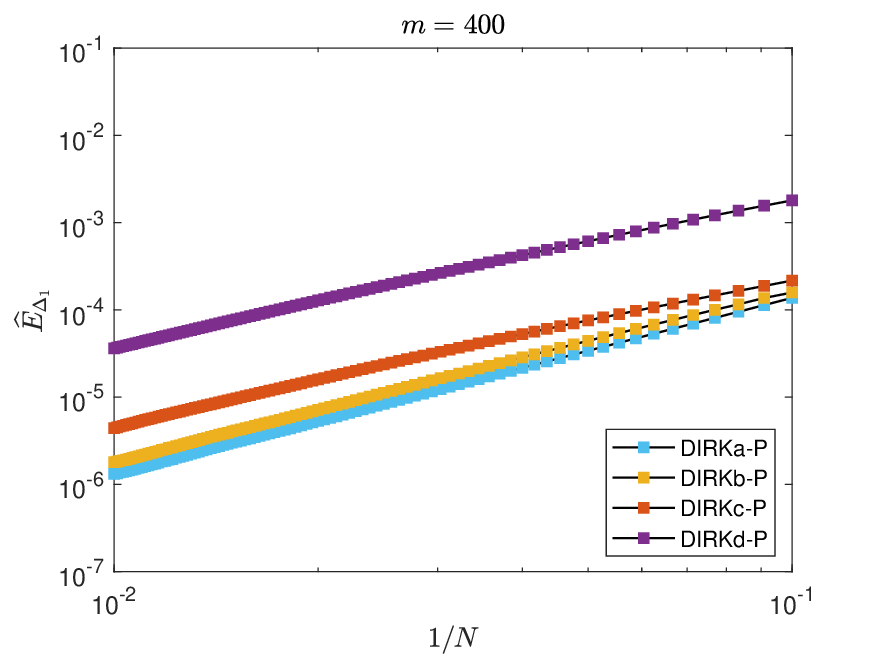}
	\includegraphics[scale=0.5]{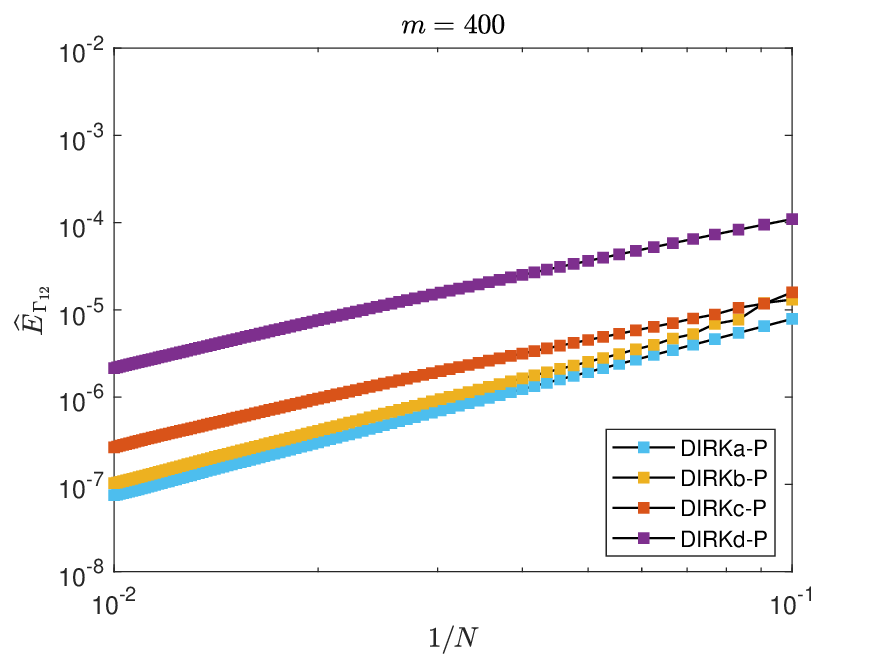}\\
	\includegraphics[scale=0.5]{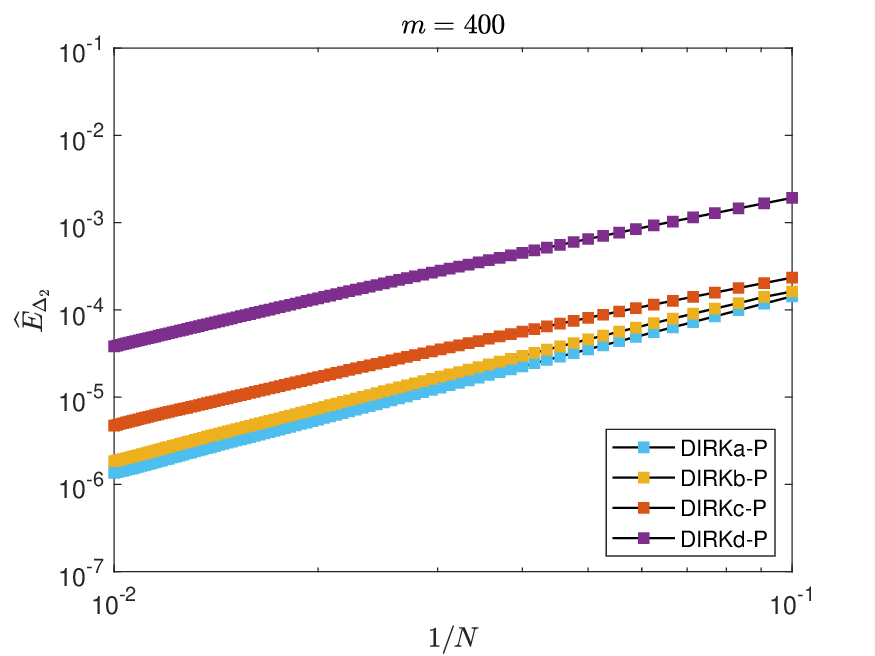}
	\includegraphics[scale=0.5]{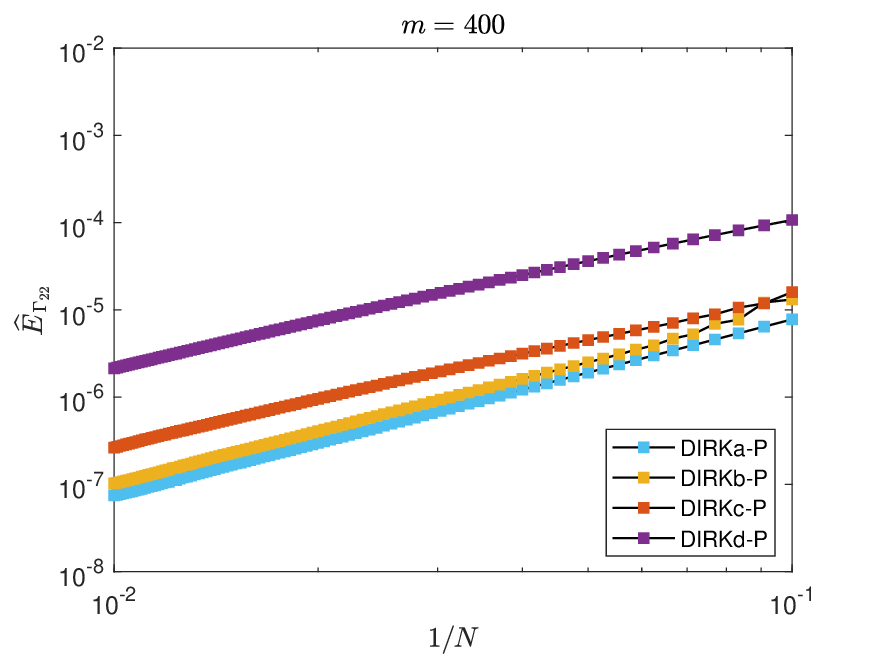}\\
	\caption{ Temporal discretization errors of the DIRKa-P, DIRKb-P, DIRKc-P, DIRKd-P methods 
	for $m=400$.
	Option value (top left), $\Delta_1$ (middle left), $\Delta_2$ (bottom left), $\Gamma_{11}$ 
	(top right), $\Gamma_{12}$ (middle right), $\Gamma_{22}$ (bottom right).}	
	\label{2Derror4}
\end{figure}

Figures~\ref{2Derror1}, \ref{2Derror2}, \ref{2Derror4} display, in double logarithmic scale,
the temporal discretization errors versus $1/N$ for $m = 100, 200, 400$, respectively, and 
$N=10,11,12,\ldots,100$.
The obtained results for the DIRKa-P, DIRKb-P, DIRKc-P and DIRKd-P methods are indicated by blue, 
orange, red and purple squares, respectively.
In each figure, the left column shows the errors in the case of the option value (top), $\Delta_1$ 
(middle) and $\Delta_2$ (bottom) and the right column shows the errors in the case of $\Gamma_{11}$ 
(top), $\Gamma_{12}$ (middle) and $\Gamma_{22}$(bottom).

A perusal of Figures~\ref{2Derror1}, \ref{2Derror2}, \ref{2Derror4} yields the positive
result that all four DIRK-P methods under consideration possess a regular, second-order temporal 
convergence behaviour for the option value as well as all Greeks Delta and Gamma, with error 
constants that are essentially independent of~$m$, which is as desired.
The error constant is found to be smallest for DIRKa-P and largest for DIRKd-P, with more than 
an order of magnitude difference.
The error constant for DIRKb-P is close to that of DIRKa-P.

We next consider the DIRKa-P method and display in Tables~\ref{tab1}, \ref{tab2}, \ref{tab3}
for $m = 100, 200, 400$, respectively, and $N=m/2$ the obtained\footnote{Spline 
interpolation is performed using neighbouring spatial grid points.} approximations for the 
option value and all Deltas and Gammas for $t=T$ at five points $(s_1,s_2)$ of interest
with $s_1, s_2 \in \{0.9K, K, 1.1K\}$. 
Computing, for each quantity and each point, the relative change in the difference between 
the approximations in Tables~\ref{tab1}, \ref{tab2} and in Tables~\ref{tab2}, \ref{tab3}, 
directly leads to the numerical convergence orders given in Table~\ref{order}.
Clearly, a favourable second-order convergence behaviour is observed for the option value 
as well as all Deltas and Gammas when $m$ and $N$ increase simultaneously in a directly 
proportional way.
For completeness, we mention that the numbers of iterations are modest; it always holds that 
$\kappa_1,\kappa_2 \in \{2,3\}$. 

\section{Conclusions}\label{Sec_Conc}
In this paper, an efficient numerical method has been derived for the two-dimensional Kou PIDCP 
to approximate the fair values of American two-asset options together with their Greeks Delta 
and Gamma under the two-asset Kou jump-diffusion model.
Here the nonlocal double integral is handled by the two-dimensional extension~\cite{HL23} of 
the fast algorithm by Toivanen~\cite{T08} for the case of the one-asset Kou model.
A useful family of DIRK methods, containing a free parameter $\theta$, has been studied for 
the temporal discretization, applied with variable step sizes.
The adaptation to the semidiscrete two-dimensional Kou PIDCP is obtained by employing the 
effective combination of the penalty iteration and the fixed-point iteration introduced by 
d'Halluin et al.~\cite{HFL04,HFV05}.

Four interesting values of $\theta$ have been considered. Numerical experiments reveal that
for each value a favourable second-order temporal convergence behaviour is obtained. 
Among the four values, we deem $\theta=1-\tfrac{1}{2}\sqrt{2}$ to be preferable in view 
of a strong stability property ($L$-stability) and a relatively small error constant of
the method.
The value $\theta=\tfrac{1}{3}$ is found to be a good alternative, provided the method 
is used with backward Euler damping.

An aim of future research is to investigate the numerical valuation of two-asset options
under infinite-activity exponential L{\'e}vy processes. 
Here among others the ideas by Wang, Wan \& Forsyth~\cite{WWF07} are expected to be fruitful.

\newpage
\clearpage

\vspace*{1cm}
\newpage

\begin{table}[h!]
\centering
\begin{tabular}{ c | c c c c c c}
 & $s_1 = 90$ & $s_1 = 100$ & $s_1 = 100$ & $s_1 = 100$ & $s_1 = 110$\\
 & $s_2 = 90$ & $s_2 = \phantom{1}90$  & $s_2 = 100$ & $s_2 = 110$ & $s_2 = 110$\\
 \hline
$u$           & ~1.4410941$\times 10^{+1}$  &  ~1.1383073$\times 10^{+1}$  &  ~8.9582601$\times 10^{+0}$  &  ~6.9719540$\times 10^{+0}$  &  ~5.2348763$\times 10^{+0}$ \\
$\Delta_1$    & -3.2584818$\times 10^{-1}$  &  -2.7943724$\times 10^{-1}$  &  -2.3504777$\times 10^{-1}$  &  -1.9394994$\times 10^{-1}$  &  -1.5432873$\times 10^{-1}$ \\
$\Delta_2$    & -3.1100634$\times 10^{-1}$  &  -2.6572941$\times 10^{-1}$  &  -2.1988068$\times 10^{-1}$  &  -1.7831686$\times 10^{-1}$  &  -1.4158342$\times 10^{-1}$ \\
$\Gamma_{11}$ & ~4.5417553$\times 10^{-3}$  &  ~4.6953152$\times 10^{-3}$  &  ~4.5270544$\times 10^{-3}$  &  ~4.1815685$\times 10^{-3}$  &  ~3.7223616$\times 10^{-3}$ \\
$\Gamma_{12}$ & ~4.4658641$\times 10^{-3}$  &  ~4.5400299$\times 10^{-3}$  &  ~4.2998511$\times 10^{-3}$  &  ~3.8964871$\times 10^{-3}$  &  ~3.4317312$\times 10^{-3}$ \\
$\Gamma_{22}$ & ~4.7484932$\times 10^{-3}$  &  ~4.7411556$\times 10^{-3}$  &  ~4.3970553$\times 10^{-3}$  &  ~3.8989366$\times 10^{-3}$  &  ~3.3963788$\times 10^{-3}$
\end{tabular}
\caption{Approximations of option value, Deltas, Gammas at five points $(s_1, s_2)$ when $m=100$, $N=50$.}
\label{tab1}
\end{table}

\begin{table}[h!]
\centering
\begin{tabular}{ c | c c c c c c}
 & $s_1 = 90$ & $s_1 = 100$ & $s_1 = 100$ & $s_1 = 100$ & $s_1 = 110$\\
 & $s_2 = 90$ & $s_2 = \phantom{1}90$  & $s_2 = 100$ & $s_2 = 110$ & $s_2 = 110$\\
 \hline
$u$           & ~1.4410326$\times 10^{+1}$  &  ~1.1382366$\times 10^{+1}$  &  ~8.9573323$\times 10^{+0}$  &  ~6.9707384$\times 10^{+0}$  &  ~5.2333520$\times 10^{+0}$ \\
$\Delta_1$    & -3.2587509$\times 10^{-1}$  &  -2.7945346$\times 10^{-1}$  &  -2.3505572$\times 10^{-1}$  &  -1.9394932$\times 10^{-1}$  &  -1.5431874$\times 10^{-1}$ \\
$\Delta_2$    & -3.1101373$\times 10^{-1}$  &  -2.6572752$\times 10^{-1}$  &  -2.1987284$\times 10^{-1}$  &  -1.7830391$\times 10^{-1}$  &  -1.4156498$\times 10^{-1}$ \\
$\Gamma_{11}$ & ~4.5418632$\times 10^{-3}$  &  ~4.6951470$\times 10^{-3}$  &  ~4.5270124$\times 10^{-3}$  &  ~4.1818159$\times 10^{-3}$  &  ~3.7229667$\times 10^{-3}$ \\
$\Gamma_{12}$ & ~4.4690948$\times 10^{-3}$  &  ~4.5427646$\times 10^{-3}$  &  ~4.3019622$\times 10^{-3}$  &  ~3.8980510$\times 10^{-3}$  &  ~3.4328213$\times 10^{-3}$ \\
$\Gamma_{22}$ & ~4.7484751$\times 10^{-3}$  &  ~4.7408826$\times 10^{-3}$  &  ~4.3969763$\times 10^{-3}$  &  ~3.8991955$\times 10^{-3}$  &  ~3.3969272$\times 10^{-3}$
\end{tabular}
\caption{Approximations of option value, Deltas, Gammas at five points $(s_1, s_2)$ when $m=200$, $N=100$.}
\label{tab2}
\end{table}

\begin{table}[h!]
\centering
\begin{tabular}{ c | c c c c c c}
 & $s_1 = 90$ & $s_1 = 100$ & $s_1 = 100$ & $s_1 = 100$ & $s_1 = 110$\\
 & $s_2 = 90$ & $s_2 = \phantom{1}90$  & $s_2 = 100$ & $s_2 = 110$ & $s_2 = 110$\\
 \hline
$u$           & ~1.4410173$\times 10^{+1}$  &  ~1.1382189$\times 10^{+1}$  &  ~8.9571007$\times 10^{+0}$  &  ~6.9704348$\times 10^{+0}$  &  ~5.2329710$\times 10^{+0}$ \\
$\Delta_1$    & -3.2588183$\times 10^{-1}$  &  -2.7945753$\times 10^{-1}$  &  -2.3505774$\times 10^{-1}$  &  -1.9394920$\times 10^{-1}$  &  -1.5431629$\times 10^{-1}$ \\
$\Delta_2$    & -3.1101559$\times 10^{-1}$  &  -2.6572706$\times 10^{-1}$  &  -2.1987090$\times 10^{-1}$  &  -1.7830070$\times 10^{-1}$  &  -1.4156039$\times 10^{-1}$ \\
$\Gamma_{11}$ & ~4.5418893$\times 10^{-3}$  &  ~4.6951043$\times 10^{-3}$  &  ~4.5270008$\times 10^{-3}$  &  ~4.1818764$\times 10^{-3}$  &  ~3.7231174$\times 10^{-3}$ \\
$\Gamma_{12}$ & ~4.4699020$\times 10^{-3}$  &  ~4.5434472$\times 10^{-3}$  &  ~4.3024889$\times 10^{-3}$  &  ~3.8984413$\times 10^{-3}$  &  ~3.4330936$\times 10^{-3}$ \\
$\Gamma_{22}$ & ~4.7484710$\times 10^{-3}$  &  ~4.7408145$\times 10^{-3}$  &  ~4.3969563$\times 10^{-3}$  &  ~3.8992599$\times 10^{-3}$  &  ~3.3970642$\times 10^{-3}$
\end{tabular}
\caption{Approximations of option value, Deltas, Gammas at five points $(s_1, s_2)$ when $m=400$, $N=200$.}
\label{tab3}
\end{table}

\begin{table}[h!]
\centering
\begin{tabular}{ c | c c c c c c}
 & $s_1 = 90$ & $s_1 = 100$ & $s_1 = 100$ & $s_1 = 100$ & $s_1 = 110$\\
 & $s_2 = 90$ & $s_2 = \phantom{1}90$  & $s_2 = 100$ & $s_2 = 110$ & $s_2 = 110$\\
 \hline
$u$           & 2.0  &  2.0  &  2.0  &  2.0  &  2.0 \\
$\Delta_1$    & 2.0  &  2.0  &  2.0  &  2.3  &  2.0 \\
$\Delta_2$    & 2.0  &  2.0  &  2.0  &  2.0  &  2.0 \\
$\Gamma_{11}$ & 2.0  &  2.0  &  1.9  &  2.0  &  2.0 \\
$\Gamma_{12}$ & 2.0  &  2.0  &  2.0  &  2.0  &  2.0 \\
$\Gamma_{22}$ & 2.1  &  2.0  &  2.0  &  2.0  &  2.0
\end{tabular}
\caption{Numerical convergence orders for option value, Deltas, Gammas at five points $(s_1, s_2)$.}
\label{order}
\end{table}

\clearpage
\bibliographystyle{plain}
\bibliography{Kou2D_AM}

\end{document}